\documentclass[11pt,reqno]{amsart}
\makeatletter
\renewcommand{\@biblabel}[1]{\makebox[4em][r]{[#1]}}
\makeatother

\usepackage{amssymb}
\usepackage{amsthm,amsmath,amssymb,amsfonts,amscd,xcolor}
\usepackage{amsfonts}
\usepackage{amsmath,amscd}
\usepackage{mathrsfs}
\usepackage{verbatim}
\usepackage{enumitem}

\usepackage{mathtools}
\usepackage[all]{xy}
\usepackage{hyperref}

\newtheorem{theorem}{Theorem}[section]
\newtheorem{lemma}[theorem]{Lemma}     
\newtheorem{corollary}[theorem]{Corollary}
\newtheorem{proposition}[theorem]{Proposition}
\newtheorem*{theorem*}{Theorem}
\newtheorem*{corollary*}{Corollary}
\newtheorem*{conjecture}{Conjecture D}

\theoremstyle{definition}

\newtheorem{remark}[theorem]{Remark}
\newtheorem{example}[theorem]{Example}

\newcommand{\mbk}{\mbox{$\Bbbk$}}

\newcommand{\mbc}{\mbox{$\mathbb{C}$}}
\newcommand{\mbq}{\mbox{$\mathbb{Q}$}}
\newcommand{\mbr}{\mbox{$\mathbb{R}$}}

\newcommand{\mfm}{\mbox{$\mathfrak{m}$}}
\newcommand{\mfn}{\mbox{$\mathfrak{n}$}}

\newcommand{\card}{\mbox{${\rm card}$}}
\newcommand{\rad}{\mbox{${\rm rad}$}}
\newcommand{\im}{\mbox{${\rm Im}$}}
\newcommand{\irr}{\mbox{${\rm Irr}$}}
\newcommand{\Spec}{\mbox{${\rm Spec}$}}

\title[Base Ideal and base locus of a Keller map]{Around the base ideal and the base locus of a Keller map}

\author[B. Enrich Lucas]{Bartomeu Enrich Lucas}
\author[F. Planas-Vilanova]{Francesc Planas-Vilanova}
\address{Departament de Matem\`atiques, Universitat Polit\`ecnica de Catalunya. 
Diagonal 647, ETSEIB, 08028 Barcelona. Catalunya}
\email{bartomeu.enrich@estudiantat.upc.edu}
\email{francesc.planas@upc.edu}
\thanks{The second author is partially supported by project PID2023-146936NB-I00, funded by the Spanish State Agency MICIU/AEI/10.13039/501100011033 and by the Generalitat de Catalunya, Departament de Recerca i Universitats, through AGAUR grant 2025 SGR 00878 (GEOMVAP)}

\date{\today}

\subjclass[2020]{14R15, 14R10, 13F20}

\begin{document}

\begin{abstract}
Let $\varphi_F:\mbk^n\to\mbk^n$ be a Keller map defined by $n$ polynomials $F_1,\ldots,F_n$ in $n$ indeterminates $X=X_1,\ldots,X_n$ over a field $\mbk$ of characteristic zero. 
Miranda-Neto defined the base ideal of $\varphi_F$ as the ideal $I_F$ generated by $F_1,\ldots,F_n$ and the base locus of $\varphi_F$ as the affine algebraic set $V_F=V(I_F)$ determined by the ideal $I_F$. He proved that $I_F$ is a radical ideal, which is a finite intersection of $v_F$ maximal ideals, where
$\card(V_F)\leq v_F$, understanding that $v_F=0$, when $I_F$ is improper. 
In this paper, we first prove that $v_F\leq \dim_{\Bbbk}(\mbk[X]/I_F)\leq [\mbk(X):\mbk(F)]$.
Using invertible affine transformations on the recent counterexample to the Jacobian Conjecture announced by Alp\"oge, we show that $I_F$ might indeed be improper 
and that the following cases, $0<\card(V_F)=v_F<[\mbk(X):\mbk(F)]$, $0<\card(V_F)=v_F=[\mbk(X):\mbk(F)]$, $0<\card(V_F)<v_F<[\mbk(X):\mbk(F)]$ and $0=\card(V_F)<v_F<[\mbk(X):\mbk(F)]$
might occur. We finish by giving a purely topological equivalent reformulation of the Jacobian Conjecture in terms of the irreducibility of $V_F$.  
\end{abstract}

\maketitle 

\section{Introduction}

Let $\mbk$ be a field of characteristic zero, $X=X_1,\ldots,X_n$ be $n$ variables over $\mbk$ and $\mbk[X]$ be the polynomial ring. Let $F_1,\ldots,F_n$ be $n$ polynomials of $\mbk[X]$ and $F=(F_1,\ldots,F_n)\in\mbk[X]^n$. Let $\varphi_F:\mbk^n\to\mbk^n$ be the polynomial map defined by $\varphi_F(x)=(F_1(x),\ldots,F_n(x))$, where $x\in\mbk^n$. The map 
$\varphi_F$ is said to be a \emph{Keller map} if the Jacobian matrix $J_F=(\partial F_i/\partial X_j)_{i,j}$ has nonzero constant determinant.
If $\varphi_F$ is an automorphism of $\mbk^n$, i.e., if $\varphi_F$ has a polynomial inverse, then, by the chain rule, $\det(J_F)$ is nonzero constant, so $\varphi_F$ is Keller. The Jacobian Conjecture states that if $\varphi_F$ is Keller, then $\varphi_F$ is an automorphism of $\mbk^n$ (see
\cite{keller}, \cite{bcw} or \cite{vde1}, for instance; 
see also \cite{rd}, for a very recent contribution to the origin of the conjecture). 
Clearly, the conjecture is true when $n=1$, for affinities or for de Jonqui\`eres maps (see, e.g., \cite{vde1}). 

Related to Keller maps and to the Jacobian Conjecture, in \cite{mn}, Miranda-Neto defines 
the \emph{base ideal $I_F$ of a Keller map} $\varphi_F$ 
as the ideal $I_F=\langle F_1,\ldots,F_n\rangle$ of $\mbk[X]$. He also defines the 
\emph{base locus $V_F$ of} $\varphi_F$ as the corresponding algebraic set $V_F=V(I_F)=V(F_1,\ldots,F_n)=\varphi_F^{-1}(\{0\})$. 
His main result, \cite[Theorem~3.1]{mn}, says that either $I_F=\mbk[X]$, or else, 
$I_F$ is a radical ideal which is the 
intersection of a finite number $v_F$ of (distinct) maximal ideals,
where $\card(V_F)\leq v_F$ and $\card(V_F)\leq [\mbk(X):\mbk(F)]$ (see
the comment before Corollary~3.2 and Remark~3.3 in \cite{mn} and 
\cite[Theorem~1.1.32]{vde1}).
In particular, if $\mbk$ is algebraically closed, then either $V_F=\varphi_F^{-1}(\{0\})=\emptyset$, which is equivalent to $I_F=\mbk[X]$, or else $V_F=\varphi_F^{-1}(\{0\})=\{a_1,\ldots,a_{v_F}\}$ is a finite set, with  $a_i=(a_{i,1},\ldots,a_{i,n})\in\mbk^n$, which is equivalent to $I_F=\cap_{i=1}^{v_F}\langle X_1-a_{i,1},\ldots,X_n-a_{i,n}\rangle$. 
Motivated by his main result, 
he states the following equivalent re-formulation of the Jacobian Conjecture:
``If $\varphi_F$ is a Keller map such that $I_F\subseteq \langle X_1,\ldots,X_n\rangle$, then 
$I_F=\langle X_1,\ldots,X_n\rangle$'' (see \cite[Conjecture C and Theorem~5.1]{mn}). 

This last summer Levent Alp\"oge \cite{alpoge} announced a counterexample to the conjecture, found with the help of AI tools, in full generality for $n=3$ and, hence, for $n\geq 4$. In \cite{ai}, there is an explanation of the reason of the failure of the conjecture. In a very recent work of Arno van den Essen in \cite{vde2} there is a very beautiful construction of a new counterexample to the Jacobian Conjecture for $n=3$. (So far, to our knowledge, the Jacobian Conjecture is still open for $n=2$.)
The example given by Alp\"oge is the following (see \url{https://x.com/jdlichtman/status/2079066717762863249}): 

\begin{example}\label{example:ccj}
Let $\mbk[X,Y,Z]$ be the polynomial ring in the variables $X,Y,Z$ over a field $\mbk$ of characteristic zero. Let $F=(F_1,F_2,F_3)$ in $\mbk[X,Y,Z]^3$, defined by
\begin{eqnarray*}
&&F_1=(1+XY)^3Z+Y^2(1+XY)(4+3XY);\\
&&F_2=Y+3X(1+XY)^2Z+3XY^2(4+3XY);\\ %F_2=Y(1-2XY)+3X(1+XY)^2Z+3XY^2(4+3XY)
&&F_3=2X-3X^2Y-X^3Z.
\end{eqnarray*}
Let $\varphi_F:\mbk^3\to\mbk^3$ be defined by $\varphi_F(x,y,z)=(F_1(x,y,z),F_2(x,y,z),F_3(x,y,z))$. 
One can check that
$\varphi_F(0,0,-1/4)=\varphi_F(1, -3/2, 13/2)=\varphi_F(-1, 3/2, 13/2)=(-1/4,0,0)$.
So, $\varphi_F$ is not injective. Moreover, one can check that $\det(J_F)=-2$, so $\varphi_F$ is a Keller map.
\end{example}

\begin{remark}
The Computer Algebra System {\sc Singular} \cite{dgps}, 
which can be used online in \url{https://singular-in-browser.pages.dev/}, allows us to check Example~\ref{example:ccj} easily. 
Namely,
\begin{verbatim}
ring A=0,(x,y,z),dp; //define the ring k[x,y,z]
poly f1=(1+x*y)^3*z+y^2*(1+x*y)*(4+3*x*y); //define f1
poly f2=y+3*x*(1+x*y)^2*z+3*x*y^2*(4+3*x*y); //define f2
poly f3=2*x-3*x^2*y-x^3*z; //define f3
ideal I=f1,f2,f3;
subst(I,x,0,y,0,z,-1/4); // image of point (0,0,-1/4)
subst(I,x,1,y,-3/2,z,13/2); // image of point (1,-3/2,13/2)
subst(I,x,-1,y,3/2,z,13/2); // image of point (-1,3/2,13/2)
matrix J=jacob(I); det(J); //Jacobian matrix and its determinant
\end{verbatim}
\end{remark}

This short note is motivated by the work of Miranda-Neto in \cite{mn}, dated in 2019, hence, prior 
to the discovery that the Jacobian Conjecture is false in dimension $n\geq 3$.  
Among other things, it caught our attention his comment
before Corollary~3.2 of \cite{mn}, where one can read 
``\emph{we set $v_F=0$, if $I_F=\mbk[X]$, even though this situation would be uninteresting and we do not believe it to be possible in the Keller case}''.

We first prove in Theorem~\ref{theo:mpd} that $v_F\leq \dim_{\Bbbk}(\mbk[X]/I_F)\leq
[\mbk(X):\mbk(F)]$. As a consequence, it follows that if $v_F=[\mbk(X,Y,Z):\mbk(F)]$, 
then $\card(V_F)=v_F$.

We point out in Corollary~\ref{cor:1=v=V} that the Keller map of Example~\ref{example:ccj} satisfies $\card(V_F)=v_F=1$, whereas $[\mbk(X,Y,Z):\mbk(F)]=3$. 

Subsequently, applying some invertible affine transformations on Example~\ref{example:ccj}, we prove in Proposition~\ref{prop:IFnoPropi} that the base ideal of a Keller map $G$ may be, in fact, the whole ring. In this case, $\card(V_G)=v_G=0$, whereas $[\mbk(X,Y,Z):\mbk(G)]=3$.
In Proposition~\ref{prop:IF<xyz} we give a counterexample to \cite[Conjecture C]{mn}, for $n=3$, by finding a non injective Keller map $\varphi_G$ with $I_G\subsetneq \langle X,Y,Z\rangle$ and 
$\card(V_G)=v_G=[\mbk(X):\mbk(G)]=3$. 
In Proposition~\ref{prop:2max}, we prove that
if $\mbk=\mbr$, then there exists a non injective Keller map $\varphi_G$ such that 
$0<\card(V_G)<v_G<[\mbk(X):\mbk(G)]=3$. In Proposition~\ref{prop:casQ}, we show that if $\mbk=\mbq$, then there exists a non surjective Keller map $\varphi_G$ such that $0=\card(V_G)<v_G=1<[\mbk(X):\mbk(G)]=3$.
Summarizing:
\vspace*{0.2cm}
\begin{center}
\begin{tabular}{c|c|c|c|l}
$\card(V_G)$   & $v_G$ & $[\mbk(X,Y,Z):\mbk(G)]$ && \\\hline
1 & 1 & 3 & Corollary~\ref{cor:1=v=V}&\mbox{ any field of characteristic zero}\\
0 & 0 & 3 & Proposition~\ref{prop:IFnoPropi}&\mbox{ any field of characteristic zero}\\
3 & 3 & 3 & Proposition~\ref{prop:IF<xyz}&\mbox{ any field of characteristic zero}\\
\hline
1 & 2 & 3 & Proposition~\ref{prop:2max} & \; $\mbk=\mbr$\\ \hline
0 & 1 & 3 & Proposition~\ref{prop:casQ} & \; $\mbk=\mbq$. 
\end{tabular}
\end{center}
\vspace*{0.2cm}

We finish the paper by giving in Theorem~\ref{theo:CD=CJ} a purely topological equivalent reformulation 
of the Jacobian Conjecture in terms of the irreducibility of the base locus. 

\section{Preliminaries}

We start with some preliminary results, some of them well-known. We give some proofs for the sake of completeness. Similar to $X$, we denote by $Y=Y_1,\ldots,Y_n$, a set of $n$ variables over $\mbk$.

\begin{remark}\label{rem:equiv}
An element $F=(F_1,\ldots,F_n)\in\mbk[X]^n$ defines at the same time a polynomial map $\varphi_F:\mbk^n\to\mbk^n$ and a morphism of $\mbk$-algebras $\phi_F:\mbk[X]\to\mbk[X]$ uniquely determined by $\phi_F(X_i)=F_i$. Let
$G=(G_1,\ldots,G_n)\in\mbk[X]^n$. Clearly,  
$\varphi_{G}\circ\varphi_F={\rm id}_{\Bbbk^n}\Leftrightarrow G(F)=(X)\Leftrightarrow 
\phi_G\circ\phi_F={\rm id}_{\Bbbk[X]}$. Thus, $\varphi_F$ is an automorphism of
$\mbk^n$ if and only if $\phi_F$ is an automorphism of $\mbk[X]$. Let us see that this is also equivalent to $\mbk[F]=\mbk[X]$.

Suppose that $\varphi_F$ is an automorphism. Then, $\varphi_G\circ\varphi_F={\rm id}_{\Bbbk^n}=
\varphi_F\circ\varphi_G$, for some $G\in\mbk[Y]^n$. Thus, $G(F)=(X)$ and $F(G)=(Y)$. In particular, 
$X_i=G_i(F_1,\ldots,F_n)\in\mbk[F_1,\ldots,F_n]$ and 
$\mbk[X]\subseteq\mbk[F]$. By definition, $\mbk[F]\subseteq \mbk[X]$. Therefore, $\mbk[F]=\mbk[X]$.
Conversely, suppose that $\mbk[F]=\mbk[X]$. Then, $F_1,\ldots,F_n$ are algebraically independent over $\mbk$. Indeed, if $H\in\mbk[X_1,\ldots,X_n]$, $H\neq 0$, is such that 
$H(F_1,\ldots,F_n)=0$, then $\phi_F(H)=H(F)=0$ and $H\in\ker(\phi_F)$. Hence, $\mbk[X]/\ker(\phi_F)\cong\im(\phi_F)=\mbk[F]$ and
$\dim(\mbk[F])=\dim(\mbk[X]/\ker(\phi))<\dim(\mbk[X])$, a contradiction. Since $\mbk[F]=\mbk[X]$, it follows 
that, for each $i=1,\ldots,n$, 
there exists $G_i\in\mbk[Y]$ such that $X_i=G_i(F_1,\ldots,F_n)=G_i(F)$. So, $G(F)=(X)$ and $\varphi_G\circ\varphi_F={\rm id}_{\Bbbk^n}$. Let us see $F(G(Y))=(Y)$, so that $\varphi_F\circ\varphi_G={\rm id}_{\Bbbk^n}$. For each
$i=1,\ldots,n$, let $H_i(Y_1,\ldots,Y_n)=F_i(G_1(Y),\ldots,G_n(Y))-Y_i\in\mbk[Y]$. 
Evaluating $Y_1,\ldots,Y_n$ in $F_1,\ldots,F_n$ and using that $G_i(F)=X_i$, we obtain
\begin{eqnarray*}
H_i(F_1,\ldots,F_n)=F_i(G_1(F),\ldots,G_n(F))-F_i=F_i(X_1,\ldots,X_n)-F_i(X_1,\ldots,X_n)=0.
\end{eqnarray*}
Since $F_1,\ldots,F_n$ are algebraically independent and $H_i(F)=0$, then
$H_i=0$ and $F(G(Y))=(Y)$.
\end{remark}

\begin{lemma}\label{lem:IFmax}
Let $\varphi_F:\mbk^n\to\mbk^n$ be a polynomial map, $F=(F_1,\ldots,F_n)\in\mbk[X]^n$. 
\begin{itemize}
\item[$(1)$] If $V_F$ is not empty, then $I_F$ is a proper ideal of $\mbk[X]$.
\item[$(2)$] If $I_F$ is a proper ideal of $\mbk[X]$ and $\mbk$ is algebraically closed, then $V_F$ is not empty. 
\item[$(3)$] If $\varphi_F:\mbk^n\to\mbk^n$ is an automorphism of $\mbk^n$, then
$\varphi_F$ is a Keller map and there exists a unique 
$a=(a_1,\dots,a_n)\in\mbk^n$ such that $I_F=\langle X_1-a_1,\ldots,X_n-a_n\rangle$. 
In particular, $V_F=\{a\}$ is irreducible. 
\item[$(4)$] If $\varphi_F,\varphi_G:\mbk^n\to\mbk^n$ are Keller, then 
$\varphi_G\circ\varphi_F:\mbk^n\to\mbk^n$ is Keller. 
\end{itemize}
\end{lemma}
\begin{proof}
If $a\in V_F$, then $F_i(a)=0$, for each $i=1,\ldots,n$. So, $F_i\in I(\{a\})$ and
$I_F=\langle F_1,\ldots,F_n\rangle \subseteq I(\{a\})=
\langle X_1-a_1,\ldots,X_n-a_n\rangle\subsetneq\mbk[X]$. On the other hand, 
if $I_F$ is proper, then there exists a maximal ideal of $\mbk[X]$ such that 
$I_F\subseteq \mfm$. If $\mbk$ is algebraically closed, by the weak Nullstellensatz, 
there is an $a\in\mbk^n$ such that
$\mfm=\langle X_1-a_1,\ldots,X_n-a_n\rangle=I(\{a\})$. 
So, $a\in V(I(\{a\}))=V(\mfm)\subseteq V(I_F)=V_F$ and $V_F$ is not empty. 

Suppose that
$\varphi_F$ is an automorphism, in particular, by the chain rule, a Keller map. 
Since $\varphi_F:\mbk^n\to\mbk^n$ is bijective, there exists a unique $a\in\mbk^n$ such that
$\varphi_F(a)=(0,\ldots,0)$. Therefore, 
\begin{eqnarray*}
V_F=\{x\in\mbk^n\mid F_1(x)=0,\ldots,F_n(x)=0\}=\{x\in\mbk^n\mid \varphi_F(x)=0\}=
\varphi_F^{-1}(\{0\})=\{a\}.
\end{eqnarray*}
Thus, $I_{F}=\langle F_1,\ldots , F_n\rangle\subseteq I(V(F_1,\ldots,F_n))=I(V_F)=
I(\{a\})=\langle X_1-a_1,\dots,X_n-a_n\rangle$. Let us see the other containment: 
$\langle X_1-a_1,\dots,X_n-a_n\rangle\subseteq I_F$. By Remark~\ref{rem:equiv}, since
$\varphi_F$ is an automorphism, then $\mbk[F]=\mbk[X]$. Thus, for each $i=1,\ldots, n$, there exists $P_i\in\mbk[Y]$, such that
$X_i-a_i=P_i(F_1,\ldots,F_n)$. Write
$P_i=\sum_{\alpha\in\mathbb{N}^n}\mu_{i,\alpha}Y^{\alpha}=\lambda_i+H_i(Y)$, where
$\lambda_i=\mu_{i,0}\in\mbk$ and 
$H_i=\sum_{\alpha\in\mathbb{N}^n\setminus\{0\}}\mu_{i,\alpha}Y^{\alpha}\in\mbk[Y]$, where 
$H_i(0)=0$. In particular, $X_i-a_i=P_i(F)=\lambda_i+H_i(F)$ and 
$\lambda_i=-H_i(F_1,\dots,F_n)+(X_i-a_i)\in I_F+\langle X_1-a_1,\ldots,X_n-a_n\rangle=
\langle X_1-a_1,\ldots,X_n-a_n\rangle$. So, there exist polynomials
$g_{i,j}\in\mbk[X]$ such that 
$\lambda_i=\sum_{j=1}^ng_{i,j}\cdot (X_j-a_j)$. Evaluating $X_j$ at $a_j$, we get
$\lambda_i=0$. Hence, for all $i=1,\dots,n$,  $X_i-a_i=H_i(F_1,\dots,F_n)\in I_F$.
Moreover, if $I_F=I(\{a\})=I(\{b\})$, then
$\{a\}=V(I(\{a\}))=V(I(\{b\}))=\{b\}$. Finally, since $I_F=I(\{a\})$, then
$V_F=V(I_F)=V(I(\{a\}))=V(X_1-a_1,\ldots,X_n-a_n)=\{a\}$ and $V_F$ is irreducible. 

The chain rule says $J_{G(F)}(X)=J_G(F(X))\cdot J_F(X)$. The entries of the matrix
$J_G=(a_{i,j})\in{\rm M}_{n\times n}(\mbk[X])$ are polynomials
$a_{i,j}=\tfrac{\partial G_i}{\partial X_j}\in\mbk[X]$. 
By definition, if $p,q\in\mbk[X]$, then $(p+q)(F)=p(F)+q(F)$ and $(pq)(F)=p(F)q(F)$. Thus, 
\begin{eqnarray*}
\det(J_G)(F)&=&
\left(\sum_{\sigma\in {\rm S}_n}(-1)^{\epsilon(\sigma)}a_{1,\sigma(1)}\cdots a_{n,\sigma(n)}\right)(F)
\stackrel{\rm def}{=}\sum_{\sigma\in {\rm S}_n}(-1)^{\epsilon(\sigma)}\left(\left(a_{1,\sigma(1)}\cdots a_{n,\sigma(n)}\right)(F)\right)\\&\stackrel{\rm def}{=}&
\sum_{\sigma\in {\rm S}_n}(-1)^{\epsilon(\sigma)}a_{1,\sigma(1)}(F)\cdots a_{n,\sigma(n)}(F)=
\det(J_G(F)).
\end{eqnarray*}
If $\varphi_F$ and $\varphi_G$ are Keller, then $\det(J_F)=c$ and $\det(J_G)=d$, where $c,d\in\mbk\setminus\{0\}$. Therefore, $\det(J_{G(F)})=\det(J_G(F))\det(J_F)
=(\det(J_G)(F))\cdot\det(J_F)=d\cdot c\in\mbk\setminus\{0\}$ and
$\varphi_G\circ\varphi_F$ is Keller.
\end{proof}

\begin{proposition}\label{prop:algext}
Let $F=(F_1,\ldots,F_n)$, $P=(P_1,\ldots,P_n)$, $Q=(Q_1,\ldots,Q_n)$, where
$F,P,Q\in\mbk[X]^n$. Let $G=Q(F(P))$, where $G=(G_1,\ldots,G_n)\in\mbk[X]^n$. Suppose that $\varphi_F$ is a Keller map and that $\varphi_P,\varphi_Q$ are invertible affine transformations. The following hold.
\begin{itemize}
\item[$(1)$] $\mbk[Q(F)]=\mbk[F]$ and $\phi_P$ induces an automorphism of  
$\mbk(X)$ such that $\phi_P(\mbk(F))=\mbk(F(P))$. 
\item[$(2)$] $F_1,\ldots,F_n$ are algebraically independent and $\mbk(F)\subset\mbk(X)$ is a finite
field extension.
\item[$(3)$] For every $x\in\mbk^n$, $\card(\varphi_F^{-1}(\{x\}))\leq [\mbk(X):\mbk(F)]$.
\item[$(4)$] $\varphi_G=\varphi_Q\circ\varphi_F\circ\varphi_P:\mbk^n\to\mbk^n$ is a Keller map
and $[\mbk(X):\mbk(G)]=[\mbk(X):\mbk(F)]$.
\end{itemize}
\end{proposition}
\begin{proof}
By hypothesis, $Q_i=\sum_{j=1}^na_{i,j}X_j+b_i$, where $a_{i,j},b_i\in\mbk$ and $M=(a_{i,j})$ is a
$n\times n$ invertible matrix. Thus, $Q_i(F)=Q_i(F_1,\ldots,F_n)=\sum_{j=1}^na_{i,j}F_j+b_i\in\mbk[F]$ and
$\mbk[Q_1(F),\ldots,Q_n(F)]\subseteq \mbk[F]$. Let $N=(c_{i,j})=M^{-1}$.
Then, $F_i=\sum_{j=1}^nc_{i,j}(Q_j(F)-b_j)\in\mbk[Q_1(F),\ldots,Q_n(F)]=\mbk[Q(F)]$. Hence, 
$\mbk[Q(F)]=\mbk[F]$. By hypothesis, $\phi_P:\mbk[X]\to\mbk[X]$ is an automorphism, with
$\phi_P(g)=g(P)$, $g\in\mbk[X]$. Therefore, it induces an automorphism 
$\phi_P:\mbk(X)\to\mbk(X)$, defined by $\phi_P(g/h)=g(P)/h(P)$. Since 
$\phi_P(F_i)=F_i(P)$, then $\phi_P(\mbk(F))\subseteq \mbk(F(P))$. Furthermore, since any $F_i(P)=\phi_P(F_i)\in\phi_P(\mbk(F))$, it follows that $\mbk(F(P))\subseteq\phi_P(\mbk(F))$.
This shows item $(1)$. The proofs of $(2)$ and $(3)$ can be found in 
\cite[Proposition~1.1.31 and Theorem~1.1.32]{vde1}. By Remark~\ref{rem:equiv}, 
$\varphi_G$ is a Keller map. By item $(1)$, $\mbk[Q(F)]=\mbk[F]$, thus
$\mbk(Q(F))=\mbk(F)$ and $[\mbk(X):\mbk(F)]=[\mbk(X):\mbk(Q(F))]$. Applying the second part of item $(1)$ to $Q(F)$, instead of $F$, we deduce that there exists an automorphism of $\mbk(X)$ such that
$\phi_P(\mbk(Q(F)))=\mbk(Q(F(P)))=\mbk(G)$. Hence, 
$[\mbk(X):\mbk(F)]=[\mbk(X):\mbk(Q(F))]=[\mbk(X):\mbk(G)]$.
\end{proof}

Recall that, by \cite[Theorem~3.1]{mn}, the base ideal $I_F$ of a Keller map $\varphi_F:\mbk^n\to\mbk^n$
is equal to the intersection of $v_F$ distinct maximal ideals, where $v_F=0$ if $I_F=\mbk[X]$.

\begin{theorem}\label{theo:mpd}
Let $\varphi_F:\mbk^n\to\mbk^n$ be a Keller map. Then, 
\begin{eqnarray*}
\card(V_F)\leq v_F\leq \dim_{\Bbbk}(\mbk[X]/I_F)\leq [\mbk(X):\mbk(F)].
\end{eqnarray*}
If $\mbk$ is algebraically closed, then $\card(V_F)=v_F$. 
If $v_F=[\mbk(X):\mbk(F)]$, then $\card(V_F)=v_F$. 
\end{theorem}
\begin{proof}
Let $I_F=\cap_{i=1}^{v_F}\mfm_i$ be the minimal primary decomposition of $I_F$, where 
$I_F=\mbk[X]$ if $v_F=0$. Then, $V_F=V(I_F)=V(\cap_{i=1}^{v_F}\mfm_i)=\cup_{i=1}^{v_F}V(\mfm_i)$. 
Note that, if $V(\mfm_i)$ is not empty, then there is $a\in V(\mfm_i)$, so $I(\{a\})\supseteq I(V(\mfm_i))\supseteq \mfm_i$. Since $\mfm_i$ is maximal, then $\mfm_i=I(\{a\})$ and
$V(\mfm_i)=V(I(\{a\})=\{a\}$ is a point. So, $V(\mfm_i)$ is either empty, or else, a point. Therefore, $\card(V_F)\leq v_F$ (see the comment just before Corollary~3.2 in \cite{mn}). In particular, if $\mbk$ is algebraically closed, then $\card(V_F)=v_F$. 

Let $L$ be the algebraic closure of $\mbk$. Let $\varphi_{F,L}:L^n\to L^n$ be the Keller map defined by $F\in L[X]^n$. Let $I_{F,L}=\langle F_1,\ldots,F_n\rangle_{L[X]}=I_FL[X]$ be the base ideal of $\varphi_{F,L}$ which, by the theorem of Miranda-Neto \cite[Theorem~3.1]{mn} can be written as $I_{F,L}=\cap_{j=1}^{v_{F,L}}\mfn_j$, where $\mfn_j$ are maximal ideals of $L[X]$. 
Since $L$ is algebraically closed, $\mfn_j=I(\{q_j\})$, for some point $q_j\in L^n$
and $L[X]/\mfn_j\cong L$. Let $V_{F,L}=V(I_{F,L})=\cup_{j=1}^{v_{F,L}}V(\mfn_j)=\{q_1,\ldots,q_{v_{F,L}}\}$ be the base locus of $\varphi_{F,L}$. 

Since $L$ is a free (hence, flat) $\mbk$-module, tensoring with $L$ preserves exactness. Thus, from
the short exact sequence $0\to I_F\to\mbk[X]\to \mbk[X]/I_F\to 0$, we deduce the short exact sequence
\begin{eqnarray*}
0\to I_F\otimes_{\Bbbk}L\to\mbk[X]\otimes_{\Bbbk}L\to (\mbk[X]/I_F)\otimes_{\Bbbk}L\to 0. 
\end{eqnarray*}
Through the natural isomorphism $\mbk[X]\otimes_{\Bbbk}L\cong L[X]$, $I_F\otimes_{\Bbbk}L$ is identified with $I_FL[X]=I_{F,L}$. Thus, $(\mbk[X]/I_F)\otimes_{\Bbbk}L\cong L[X]/I_{F,L}$. Since
$\mbk[X]/I_F$ is a free $\mbk$-module with a basis $\{e_i\}_{i\in\Lambda}$, say, then
$L[X]/I_{F,L}\cong (\mbk[X]/I_F)\otimes_{\Bbbk}L$ is a free $L$-module with a basis
$\{e_i\otimes 1\}_{i\in\Lambda}$. In particular, 
$\dim_{\Bbbk}(\mbk[X]/I_F)=\dim_L(L[X]/I_{F,L})$.
By the Chinese Remainder Theorem, 
\begin{eqnarray*}
\mbk[X]/I_F\cong \mbk[X]/\mfm_1\times\ldots\times\mbk[X]/\mfm_{v_F}\mbox{ and }L[X]/I_{F,L}\cong L[X]/\mfn_1\times\ldots\times L[X]/\mfn_{v_{F,L}}.%\cong \prod_{j=1}^{v_{F,L}}L.
\end{eqnarray*}
In particular, 
\begin{eqnarray*}
v_F\leq \sum_{i=1}^{v_F}\dim_{\Bbbk}(\mbk[X]/\mfm_i)=\dim_{\Bbbk}(\mbk[X]/I_F)=
\dim_L(L[X]/I_{F,L})=\sum_{j=1}^{v_{F,L}}\dim_L(L[X]/\mfn_j)=v_{F,L}.
\end{eqnarray*}
Since $L$ is algebraically closed, $v_{F,L}=\card(V_{F,L})$. By Proposition~\ref{prop:algext}, $\card(V_{F,L})\leq [L(X):L(F)]$.
 
Observe that $\mbk(X)=\mbk(F)(X_1,\ldots,X_n)$ and 
$L(X)=L(F)(X_1,...,X_n)$. For each $i=1,\ldots,n$, let $M_i=\mbk(F)(X_1,\ldots,X_i)$ and $N_i=L(F)(X_1,\ldots,X_i)$, $M_0=\mbk(F)$ and $N_0=L(F)$. In particular, $M_n=\mbk(X)$ and $N_n=L(X)$. Since $X_1,\ldots,X_n$ are algebraic over $\mbk(F)$, it follows that, for each $i=1,\ldots,n$, $M_{i-1}\subset M_{i}$ and $N_{i-1}\subset N_{i}$ are algebraic extensions and $[N_{i}:N_{i-1}]=\deg(\irr(X_{i},N_{i-1}))\leq\deg(\irr(X_{i},M_{i-1}))=[M_i:M_{i-1}]$, where $\irr(\alpha,K)$ stands for the irreducible polynomial of an algebraic element $\alpha$ over a field $K$.
Thus,
\begin{eqnarray*}
[L(X):L(F)]=\prod_{i=1}^{n}[N_{i}:N_{i-1}]\leq \prod_{i=1}^{n}[M_{i}:M_{i-1}]=[\mbk(X):\mbk(F)].
\end{eqnarray*}
Therefore, 
$v_F\leq \dim_{\Bbbk}(\mbk[X]/I_F)=v_{F,L}=\card(V_{F,L})\leq [L(X):L(F)]\leq [\mbk(X):\mbk(F)]$.
In particular, if $v_F=[\mbk(X):\mbk(F)]$, then 
$v_F=\dim_{\Bbbk}(\mbk[X]/I_F)=\sum_{i=1}^{v_F}\dim_{\Bbbk}(\mbk[X]/\mfm_i)$. Thus, each
$\mbk[X]/\mfm_i$ has dimension $1$ and $\mfm_i=I(\{p_i\})$, for some point $p_i\in\mbk^n$. 
So, $V_F=V(I_F)=\cup_{i=1}^{v_F}V(\mfm_i)=\{p_1,\ldots,p_{v_{F}}\}$ and $\card(V_F)=v_F$.
\end{proof}

\section{The degree of the field extension in Example~\ref{example:ccj}}

Let us calculate the degree of the algebraic field extension $\mbk(F_1,F_2,F_3)\subset \mbk(X,Y,Z)$ in Example~\ref{example:ccj}. See \cite[page~5, before Remark~4.4]{ai}, for a somewhat different approach.

\begin{proposition}\label{prop:degree}
Let $F_1,F_2,F_3$ be the three polynomials in {\rm Example~\ref{example:ccj}} and $F=(F_1,F_2,F_3)$.
Let $A=27Y_1^2Y_3^2-18Y_1Y_2Y_3+16Y_1+Y_2^3Y_3-Y_2^2$, $B=4-3Y_2Y_3$ and $C=-2Y_3$, $A,B,C\in\mbk[Y_1,Y_2,Y_3]$. Then, 
\begin{eqnarray}
&&A(F_1,F_2,F_3)X^3+B(F_1,F_2,F_3)X+C(F_1,F_2,F_3)=0,\label{eq:x}\\
&&Y=(9X^2F_1F_3-3XF_2F_3-X^2F_2-3F_3+6X)/(3X^2F_2F_3+3XF_3-4X^2)\mbox{ and }\label{eq:y}\\
&&Z=(2X-3X^2Y-F_3)/X^3\label{eq:z}.
\end{eqnarray}
In particular, $[\mbk(X,Y,Z):\mbk(F)]=3$. 
\end{proposition}
\begin{proof}
Recall that, by definition, $F_1,F_2,F_3$ are:
$F_1=(1+XY)^3Z+Y^2(1+XY)(4+3XY)$, $F_2=Y+3X(1+XY)^2Z+3XY^2(4+3XY)$ and 
$F_3=2X-3X^2Y-X^3Z$. Thus, 
\begin{eqnarray*}%\label{eq:z}
Z=(2X-3X^2Y-F_3)/X^3. 
\end{eqnarray*}
It follows that,
$\mbk(X,Y,Z)=\mbk(X,Y,F_1,F_2,F_3)$. On the other hand, one can check that:
\begin{eqnarray}\label{eq:A}
3XF_1-(1+XY)(F_2-Y)=0.
\end{eqnarray}    
Multiplying $F_2$ by $X^2$, we obtain:
\begin{eqnarray}\label{eq:x2f2}
X^2F_2=X^2Y+3(1+XY)^2(X^3Z)+3X^3Y^2(4+3XY).
\end{eqnarray}
Substituting $X^3Z$ by $2X-3X^2Y-F_3$ in \eqref{eq:x2f2}, we get:
\begin{eqnarray*}
X^2F_2=-3X^2Y^2F_3+(4X^2-6XF_3)Y+(6X-3F_3),
\end{eqnarray*}
so 
\begin{eqnarray}\label{eq:B}
3X^2Y^2F_3-(4X^2-6XF_3)Y-(6X-X^2F_2-3F_3)=0.
\end{eqnarray}
Multiplying equation \eqref{eq:A} by $3XF_3$ and subtracting equation \eqref{eq:B}, we get:
\begin{eqnarray*}
3XF_3(3XF_1-(1+XY)(F_2-Y))-3X^2Y^2F_3+(4X^2-6XF_3)Y+(6X-X^2F_2-3F_3)=0.
\end{eqnarray*}
Developing this expression and collecting the terms in the variable $Y$:
\begin{eqnarray}\label{eq:Y}
Y(-3X^2F_2F_3-3XF_3+4X^2)+(9X^2F_1F_3-3XF_2F_3-X^2F_2-3F_3+6X)=0.
\end{eqnarray}
Let $N=9X^2F_1F_3-3XF_2F_3-X^2F_2-3F_3+6X$ and $D=3X^2F_2F_3+3XF_3-4X^2$. Then, $Y=N/D$
and $Y\in \mbk(X,F_1,F_2,F_3)$. So $\mbk(X,Y,Z)=\mbk(F_1,F_2,F_3)(X)$. From equation \eqref{eq:A}, we get:
\begin{eqnarray}\label{eq:A'}
XY^2+(1-XF_2)Y+(3XF_1-F_2)=0.
\end{eqnarray}
Substituting $Y=N/D$ in equation~\eqref{eq:A'} and multiplying by $D^2$, we get:
\begin{eqnarray}\label{eq:A'Y}
XN^2+(1-XF_2)DN+(3XF_1-F_2)D^2=0. 
\end{eqnarray}
Substituting the values of $D$ and $N$ in equation~\ref{eq:A'Y}, we get: 
\begin{eqnarray*}
3X^2\cdot\left[(27F_1^2F_3^2-18F_1F_2F_3+16F_1+F_2^3F_3-F_2^2)X^3+(4-3F_2F_3)X-2F_3\right]=0.
\end{eqnarray*}
Simplifying by $3X^2$, we obtain
\begin{eqnarray}\label{eq:Xalg}
\underbrace{(27F_1^2F_3^2-18F_1F_2F_3+16F_1+F_2^3F_3-F_2^2)}_AX^3+
\underbrace{(4-3F_2F_3)}_BX+\underbrace{(-2F_3)}_{C}=0.
\end{eqnarray}
Hence, $X$ is algebraic over $\mbk(F)$ and  
$[\mbk(F)(X):\mbk(F)]\leq 3$.
Example~\ref{example:ccj} says that $(-1/4,0,0)$ has at least
3 preimages through $\varphi_F$. 
By Proposition~\ref{prop:algext}, $[\mbk(X,Y,Z):\mbk(F)]\geq 3$.
Thus, we get the equality $[\mbk(X,Y,Z):\mbk(F)]=3$.
To check equalities~\eqref{eq:Y} and \eqref{eq:Xalg} with {\sc Singular} \cite{dgps},
\url{https://singular-in-browser.pages.dev/}:

\begin{verbatim}
ring A=(0,a,b,c),(x,y,z),dp;  //define the ring A
poly f1=(1+x*y)^3*z+y^2*(1+x*y)*(4+3*x*y); //define f1 
poly f2=y+3*x*(1+x*y)^2*z+3*x*y^2*(4+3*x*y); //define f2
poly f3=2*x-3*x^2*y-x^3*z; //define f3
y*x*(3*x*f2*f3+3*f3-4*x)-(9*x^2*f1*f3-3*x*f2*f3-f2*x^2-3*f3+6*x); 
//this is zero, so y is in k(x,f1,f2,f3)
poly ny=9*x^2*a*c-3*x*b*c-b*x^2-3*c+6*x; //numerator of y
poly dy=x*(3*x*b*c+3*c-4*x); //denominator of y
poly g=x*ny^2+(1-x*b)*dy*ny+(3*x*a-b)*dy^2; 
poly h=(27*a^2*c^2-18*a*b*c+16*a+b^3*c-b^2)*x^3+(4-3*b*c)*x-2*c; 
g-3*x^2*h; //this is zero
(27*f1^2*f3^2-18*f1*f2*f3+16*f1+f2^3*f3-f2^2)*x^3+(4-3*f2*f3)*x-2*f3; //zero
\end{verbatim}
%\vspace*{-0.7cm}
\end{proof}

\begin{corollary}\label{cor:inIG}
Let $F_1,F_2,F_3$ be the three polynomials in {\rm Example~\ref{example:ccj}} and $F=(F_1,F_2,F_3)$. 
Let $A,B,C\in\mbk[Y_1,Y_2,Y_3]$ be the polynomials defined in {\rm Proposition~\ref{prop:degree}}.
Let $(a,b,c)\in\mbk^3$, $G_1=F_1-a$, $G_2=F_2-b$, $G_3=F_3-c$ and $G=(G_1,G_2,G_3)$. 
Then, $\varphi_G:\mbk^3\to \mbk^3$ is a Keller map and 
\begin{eqnarray}\label{eq:belongIF}
A(a,b,c)X^3+B(a,b,c)X+C(a,b,c)\in I_G.
\end{eqnarray}
Moreover, if $(a,b,c)=\varphi_F(x,y,z)$, for some $(x,y,z)\in\mbk^3$, then 
\begin{eqnarray}\label{eq:resultant}
A(a,b,c)x^3+B(a,b,c)x+C(a,b,c)=0,
\end{eqnarray}
where $A(a,b,c)=27a^2c^2-18abc+16a+b^3c-b^2$, $B(a,b,c)=4-3bc$ and $C(a,b,c)=-2c$.
\end{corollary}
\begin{proof}
Clearly, $\varphi_G:\mbk^3\to\mbk^3$ is a Keller map (see, also, Proposition~\ref{prop:algext}). By the equality~\eqref{eq:x} in Proposition~\ref{prop:degree} and using the binomial expansion, we get:
\begin{eqnarray*}
0&=&A(F_1,F_2,F_3)X^3+B(F_1,F_2,F_3)X+C(F_1,F_2,F_3)\\
&=&A(G_1+a,G_2+b,G_3+c)X^3+B(G_1+a,G_2+b,G_3+c)X+C(G_1+a,G_2+b,G_3+c)\\
&=&A(a,b,c)X^3+B(a,b,c)X+C(a,b,c)+PX^3+QX+R,
\end{eqnarray*}
where $P,Q,R\in\langle G_1,G_2,G_3\rangle=I_G$, so 
$A(a,b,c)X^3+B(a,b,c)X+C(a,b,c)\in I_G$.

Suppose that $a=F_1(x,y,z)$, $b=F_2(x,y,z)$ and $c=F_3(x,y,z)$, for some $(x,y,z)\in\mbk^3$. Substituting $X,Y,Z$ by $x,y,z$ in 
$A(F_1,F_2,F_3)X^3+B(F_1,F_2,F_3)X+C(F_1,F_2,F_3)=0$, which is an equality in $\mbk[X,Y,Z]$, we get $A(a,b,c)x^3+B(a,b,c)x+C(a,b,c)=0$. Making explicit the values of
$A(a,b,c)$, $B(a,b,c)$ and $C(a,b,c)$, we get
$(27a^2c^2-18abc+16a+b^3c-b^2)x^3+(4-3bc)x-2c=0$.
\end{proof}

\begin{corollary}\label{cor:1=v=V}
Let $F_1,F_2,F_3$ be the three polynomials in {\rm Example~\ref{example:ccj}} and $F=(F_1,F_2,F_3)$. 
Then, $I_F=\langle X,Y,Z\rangle$ and $\card(V_F)=v_F=1$, whereas $[\mbk(X,Y,Z):\mbk(F)]=3$. 
\end{corollary}
\begin{proof}
By Proposition~\ref{prop:degree}, $[\mbk(X,Y,Z):\mbk(F)]=3$. Take in equation~\eqref{eq:belongIF}, $(a,b,c)=(0,0,0)$. Then, $G=(G_1,G_2,G_3)=(F_1,F_2,F_3)=F$, $I_G=I_F$ and
$A(0,0,0)X^3+B(0,0,0)X+C(0,0,0)\in I_G=I_F$. But, $A(0,0,0)=0$, $B(0,0,0)=4$ and $C(0,0,0)=0$. Therefore, 
$4X$, and hence, $X\in I_F$. Since $F_2=Y+XP$, then $Y\in \langle F_2,X\rangle\subset I_F$. 
Since $F_1=(1+XY)^3Z+YQ$, then $F_1=Z+YR$ and $Z\in \langle F_1,Y\rangle\subset I_F$. Thus, 
$\langle X,Y,Z\rangle\subseteq I_F$. The other inclusion follows by the definition of $F_1,F_2,F_3$.
In particular, $v_F=1$, $V_F=V(I_F)=\{(0,0,0)\}$ and $\card(V_F)=v_F$. 
To check it with {\sc Singular}, \cite{dgps}, add the line
\begin{verbatim}
ideal I=f1,f2,f3; lift(I,x); lift(I,y); lift(I,z);
\end{verbatim}
to the code in Example~\ref{example:ccj}. 
\end{proof}

Corollary~\ref{cor:1=v=V} shows that there exists a Keller map $\varphi_F$ which 
is not an automorphism and such that $I_{F}=\langle X,Y,Z\rangle$. Thus, 
the converse of Lemma~\ref{lem:IFmax}, $(3)$, is false in general. 

\section{The base ideal might be improper}

To see that the base ideal of a Keller map might be improper we use the following lemma. 

\begin{lemma}\label{lem:vgnobuit}
Let $\varphi_F:\mbk^n\to\mbk^n$ be a non surjective Keller map. 
Let $a=(a_1,\ldots,a_n)\in\mbk^n$, $a\not\in\im(\varphi_F)$. Let 
$G=(F_1-a_1,\ldots,F_n-a_n)\in\mbk[X]^n$. Then, 
$\varphi_G:\mbk^n\to\mbk^n$ is a non surjective Keller map and $V_G$ is empty. Moreover,
if $\mbk$ is algebraically closed, then $I_G=\mbk[X]$ is not a proper ideal.
\end{lemma}
\begin{proof}
Let $Q=(X_1-a_1,\ldots,X_n-a_n)\in\mbk[X]^n$. Then, $G=Q(F)$ and $\varphi_G=\varphi_Q\circ\varphi_F$. 
Since $\varphi_F$ and $\varphi_Q$ are Keller maps, then $\varphi_G$ is a Keller map. 
Since $\varphi_F$ is not surjective and $\varphi_Q$ is invertible, then $\varphi_G$ is 
not surjective. In addition, since $a$ is not in the image of $\varphi_F$, 
$V(I_G)=V(F_1-a_1,\ldots,F_n-a_n)=\{x\in\mbk^n\mid F_1(x)=a_1,\ldots,F_n(x)=a_n\}=\emptyset$.
If $\mbk$ is algebraically closed, then 
$\rad(I_G)=I(V(I_G))=I(\emptyset)=\mbk[X]$ and $I_G=\mbk[X]$ is not a proper ideal.  
\end{proof}

Thus, when $\mbk$ is algebraically closed, to show that the base ideal 
of a Keller map may fail to be proper, it is enough to find a non surjective one. 
We know that the Keller map of Example~\ref{example:ccj} is not injective. 
Let us see that it is neither surjective. In fact, 
the point we find with empty fiber corresponds to $t=3/2$ in
\cite[Theorem~4.2]{ai}. Our approach here is slightly different and does not use $\mbk=\mbc$.

\begin{lemma}\label{lem:varphiFnoexh}
Let $F_1,F_2,F_3$ be the three polynomials in {\rm Example~\ref{example:ccj}} and $F=(F_1,F_2,F_3)$. 
Let $A,B,C\in\mbk[Y_1,Y_2,Y_3]$ be the polynomials defined in {\rm Proposition~\ref{prop:degree}}.
Then, the point $(\tfrac{4}{27},\tfrac{4}{3},1)$ is not in the image of $\varphi_F$. Moreover,
$A(\tfrac{4}{27},\tfrac{4}{3},1)=0$, $B(\tfrac{4}{27},\tfrac{4}{3},1)=0$ and $C(\tfrac{4}{27},\tfrac{4}{3},1)=-2$. 
\end{lemma}
\begin{proof}
If $(a,b,c)=\varphi_F(x,y,z)$, for some $(x,y,z)\in\mbk^3$, then, 
by the equality~\eqref{eq:resultant} in Corollary~\ref{cor:inIG}:
\begin{eqnarray*}
A(a,b,c)x^3+B(a,b,c)x+C(a,b,c)=0,
\end{eqnarray*}
where $A(a,b,c)=(27a^2c^2-18abc+16a+b^3c-b^2)$, $B(a,b,c)=4-3bc$ and $C(a,b,c)=-2c$. Forcing
the inconsistency of this equation in $x$, we may take $c=1$ and $b=4/3$, so that
$C(\cdot,\cdot,1)=-2$ and $B(\cdot,4/3,1)=0$. 
If, in addition, one takes $a=4/27$, then $A(4/27,4/3,1)=0$. 
\end{proof}

From Lemmas~\ref{lem:varphiFnoexh} and \ref{lem:vgnobuit}, one can deduce 
that, if $\mbk$ is algebraically closed, then 
there exists a Keller map $\varphi_G:\mbk^3\to\mbk^3$
such that $I_G=\mbk[X]$. In fact, one can even avoid the assumption on
$\mbk$ being algebraically closed. Recall that in Example~\ref{example:ccj}, $\mbk$ is just 
a field of characteristic zero. 

\begin{proposition}\label{prop:IFnoPropi}
Let $F_1,F_2,F_3$ be the three polynomials in {\rm Example~\ref{example:ccj}} and $F=(F_1,F_2,F_3)$. 
Let $G_1=F_1-\tfrac{4}{27}$, $G_2=F_2-\tfrac{4}{3}$, $G_3=F_3-1$ and $G=(G_1,G_2,G_3)$. 
Then, $\varphi_G:\mbk^3\to \mbk^3$ is a non surjective Keller map and $I_G=\mbk[X]$. 
In particular, $V_G$ is not irreducible. 
\end{proposition}
\begin{proof} 
By Example~\ref{example:ccj} and Lemma~\ref{lem:varphiFnoexh}, 
$\varphi_F$ is a non surjective Keller map.
By Lemma~\ref{lem:vgnobuit}, $\varphi_G$ is a non surjective Keller map. 
By Corollary~\ref{cor:inIG} and Lemma~\ref{lem:varphiFnoexh}, $A(\tfrac{4}{27},\tfrac{4}{3},1)X^3+B(\tfrac{4}{27},\tfrac{4}{3},1)X+C(\tfrac{4}{27},\tfrac{4}{3},1)=-2\in I_G$. Therefore, $I_G=\mbk[X]$ and $V_G=V(I_G)=\emptyset$ is not irreducible. 

With {\sc Singular} \cite{dgps}, \url{https://singular-in-browser.pages.dev/}
for the online version:
\begin{verbatim}
ring A=0,(x,y,z),dp;  //define the ring
poly f1=(1+x*y)^3*z+y^2*(1+x*y)*(4+3*x*y); poly g1=f1-(4/27); //define f1 and g1
poly f2=y+3*x*(1+x*y)^2*z+3*x*y^2*(4+3*x*y); poly g2=f2-(4/3); //define f2 and g2
poly f3=2*x-3*x^2*y-x^3*z; poly g3=f3-1; //define f3 and g3
ideal Ig=g1,g2,g3;  //define the ideal J 
lift(Ig,1);  // express 1 in terms of the gi
\end{verbatim}
%\vspace*{-0.5cm}
\end{proof}

\section{The base ideal of a Keller map might be proper and not maximal}

To see that the base ideal of a Keller map might be strictly included in the maximal
ideal at the origin we use the following lemma. 

\begin{lemma}\label{lem:notinjective}
Let $\varphi_F:\mbk^n\to\mbk^n$ be a non injective Keller map. 
Let $a,b,c\in\mbk^n$, where $a=(a_1,\ldots,a_n)$, $b=(b_1,\ldots,b_n)$, $a\neq b$, 
and $c=(c_1,\ldots,c_n)$, such that $\varphi_F(a)=\varphi_F(b)=c$.
Let 
\begin{eqnarray*}
P=(X_1+a_1,\ldots,X_n+a_n),\; Q=(X_1-c_1,\ldots,X_n-c_n)\;\mbox{ and }\; G=Q(F(P)),
\end{eqnarray*}
where $F,P,Q,G\in\mbk[X]^n$. Let $\varphi_G:\mbk^n\to \mbk^n$ be the polynomial map defined by $G$. 
Then, $\varphi_G=\varphi_Q\circ\varphi_F\circ\varphi_P$ is a non injective Keller map, 
$V_G\supseteq \{0,b-a\}$ and
\begin{eqnarray*}
I_G\subseteq\langle X_1,\ldots,X_n\rangle\cap
\langle X_1-(b_1-a_1),\ldots,X_n-(b_n-a_n)\rangle.
\end{eqnarray*}
In particular, $I_G\subsetneq \langle X_1,\ldots,X_n\rangle$.
\end{lemma}
\begin{proof}
By definition of $G$, $\varphi_G=\varphi_Q\circ\varphi_F\circ\varphi_P$.
Clearly, $\varphi_P$ and $\varphi_Q$ are invertible Keller maps.
Thus, $\varphi_G$ is a non injective Keller map. Moreover, 
$\varphi_G(0)=\varphi_Q(\varphi_F(\varphi_P(0)))=\varphi_Q(\varphi_F(a))=\varphi_Q(c)=0$. 
Similarly, $\varphi_G(b-a)=
\varphi_Q(\varphi_F(\varphi_P(b-a)))=\varphi_Q(\varphi_F(b))=\varphi_Q(c)=0$. Hence, 
$\{0,b-a\}\subseteq\varphi_G^{-1}(\{0\})=V_G=V(I_G)$ and 
\begin{eqnarray*}
I(\{0\})\cap I(\{b-a\})=I(\{0,b-a\})\supseteq I(V(I_G))\supseteq I_G.
\end{eqnarray*}
Therefore, $I_G\subseteq \langle X_1,\ldots,X_n\rangle\cap 
\langle X_1-(b_1-a_1),\ldots,X_n-(b_n-a_n)\rangle$.
Since $a\neq b$, it follows that $I_G\subsetneq \langle X_1,\ldots,X_n\rangle$.
\end{proof}

By suitably modifying Example~\ref{example:ccj} and applying Lemma~\ref{lem:notinjective}, 
we achieve a counterexample to Conjecture C, for $n=3$ (\cite[Conjecture C]{mn}). 
Here, we get $\card(V_G)=v_G=[\mbk(X,Y,Z):\mbk(G)]$. 

\begin{proposition}\label{prop:IF<xyz} 
Let $F_1,F_2,F_3$ be the three polynomials in {\rm Example~\ref{example:ccj}} and
$F=(F_1,F_2,F_3)$. Let $P=(P_1,P_2,P_3)=(X,Y,Z-\tfrac{1}{4})$, 
$Q=(Q_1,Q_2,Q_3)=(X+\tfrac{1}{4},Y,Z)$ and $G=Q(F(P))$, where
$F,P,Q,G$ are in $\mbk[X,Y,Z]^3$.
Let $\varphi_G:\mbk^3\to \mbk^3$ be the polynomial map defined by $G$. 
Then, $\varphi_G=\varphi_Q\circ\varphi_F\circ\varphi_P$ is a non injective Keller map 
and $I_G\subsetneq \langle X,Y,Z\rangle$. Concretely, 
\begin{eqnarray*}
I_G=\langle X,Y,Z\rangle\cap\langle X-1,Y+\tfrac{3}{2},Z-\tfrac{27}{4}\rangle
\cap\langle X+1,Y-\tfrac{3}{2},Z-\tfrac{27}{4}\rangle,
\end{eqnarray*}
$V_G=\{(0,0,0),(1,-3/2,27/4),(-1,3/2,27/4)\}$ and
$\card(V_G)=v_G=[\mbk(X,Y,Z):\mbk(G)]=3$.
\end{proposition}
\begin{proof} 
Set $a=(0,0,-1/4)$, $b=(1,-3/2,13/2)$, $b'=(-1,3/2,13/2)$ and $c=(-1/4,0,0)$.
Example~\ref{example:ccj} tells us that $\varphi_F(a)=\varphi_F(b)=\varphi_F(b')=c$.
By Lemma~\ref{lem:notinjective}, we deduce that 
$\varphi_G=\varphi_Q\circ\varphi_F\circ\varphi_P$ is a non injective Keller map, 
$V_G\supseteq \{0,b-a,b'-a\}$ and $I_G\subseteq\mfm_1\cap\mfm_2\cap \mfm_3$, where 
$\mfm_1=\langle X,Y,Z\rangle$, $\mfm_2=\langle X-1,Y+\tfrac{3}{2},Z-\tfrac{27}{4}\rangle$
and $\mfm_3=\langle X+1,Y-\tfrac{3}{2},Z-\tfrac{27}{4}\rangle$.
Using Propositions~\ref{prop:algext}, \ref{prop:degree} and Theorem~\ref{theo:mpd}, we deduce that
$3\leq \card(V_G)\leq v_G\leq [\mbk(X):\mbk(G)]=3$. Thus, $v_G=3$ and $I_G=\mfn_1\cap\mfn_2\cap\mfn_3\subseteq \mfm_1\cap\mfm_2\cap\mfm_3\subset \mfm_i$, so
there exists $\mfn_{j(i)}\subseteq \mfm_i$. Since $\mfn_{j(i)}$ is maximal, it follows that
$\mfn_{j(i)}=\mfm_i$ and $\{\mfn_1,\mfn_2,\mfn_3\}=\{\mfm_1,\mfm_2,\mfm_3\}$. Hence, $I_G=\mfm_1\cap\mfm_2\cap \mfm_3$.
To check this, we can use {\sc Singular} \cite{dgps}, \url{https://singular-in-browser.pages.dev/}:
\begin{verbatim}
ring A=0,(x,y,z),dp;  //define the ring A
option(redSB); // forcing the unicity of the reduced GB
poly f1=(1+x*y)^3*z+y^2*(1+x*y)*(4+3*x*y); //define f1 
poly f2=y+3*x*(1+x*y)^2*z+3*x*y^2*(4+3*x*y); //define f2
poly f3=2*x-3*x^2*y-x^3*z; //define f3
poly f1p=subst(f1,z,z-1/4); //define f_1(p)
poly f2p=subst(f2,z,z-1/4); //define f_2(p)
poly f3p=subst(f3,z,z-1/4); //define f_3(p)
poly g1=f1p+1/4; poly g2=f2p; poly g3=f3p; 
ideal Ig=g1,g2,g3; //define ideal Ig
ideal m1=x,y,z; ideal m2=x-1,y+3/2,z-27/4; ideal m3=x+1,y-3/2,z-27/4;
ideal J=intersect(m1,m2,m3); std(Ig); std(J); // Ig and J are equal
\end{verbatim}
%\vspace*{-0.6cm}
\end{proof}

\section{The real case}

If $\mbk=\mbr$, then there exists a Keller map $\varphi_G$ with
$0<\card(V_G)<v_G<[\mbk(X,Y,Z):\mbk(G)]=3$. This answers, in dimension $3$, a question of Miranda-Neto in dimension $2$ (see \cite{mn}, just before Theorem~5.1). We start by finding a point in $\mbc^3$ whose preimage by $\varphi_F$ of Example~\ref{example:ccj} has two non-real points.

\begin{lemma}\label{lem:2max}
Let $F_1,F_2,F_3$ be the three polynomials in {\rm Example~\ref{example:ccj}} and $F=(F_1,F_2,F_3)$.
Then, $\varphi_F^{-1}(\{(\frac{1}{4},0,0)\})=
\{(0,0,1/4),(i,\tfrac{3}{2}i,-\tfrac{13}{2}),(-i,-\tfrac{3}{2}i,-\tfrac{13}{2})\}$.
\end{lemma}
\begin{proof}
If $(a,b,c)=\varphi_F(x,y,z)$, for some $(x,y,z)\in\mbk^3$, then
\begin{eqnarray}\label{eq:antiimatge}
\left\{\begin{array}{ccccl}
a&=&F_1(x,y,z)&=&(1+xy)^3z+y^2(1+xy)(4+3xy),\\
b&=&F_2(x,y,z)&=&y+3x(1+xy)^2z+3xy^2(4+3xy),\\ 
c&=&F_3(x,y,z)&=&2x-3x^2y-x^3z.\end{array}\right.
\end{eqnarray} 
Moreover, by the equality~\eqref{eq:resultant} in Corollary~\ref{cor:inIG}:
\begin{eqnarray*}
A(a,b,c)x^3+B(a,b,c)x+C(a,b,c)=0,
\end{eqnarray*}
where $A(a,b,c)=(27a^2c^2-18abc+16a+b^3c-b^2)$, $B(a,b,c)=4-3bc$ and $C(a,b,c)=-2c$.
Take $\mbk=\mbc$. To have a unique root $x\in\mbr$ in $A(a,b,c)x^3+B(a,b,c)x+C(a,b,c)=0$, let
$A(a,b,c)=B(a,b,c)$ and $C(a,b,c)=0$, so that the equation becomes $A(a,b,c)\cdot x(x^2+1)=0$. 
But, if $C(a,b,c)=0$, then $c=0$ and, if moreover $A(a,b,c)=B(a,b,c)$, then 
$B(a,b,c)=4$ and $A(a,b,c)=16a-b^2=4$. This holds, for instance, for $(a,b,c)=(1/4,0,0)$. In this case, $x\in\{0,i,-i\}$. If $x=0$, using equations~\eqref{eq:antiimatge}, it follows that $1/4=a=z+4y$ and $0=b=y$, $(x,y,z)=(0,0,1/4)$. Suppose that $x=\pm i$ and that $(1/4,0,0)=\varphi_F(i,y,z)=(F_1(i,y,z),F_2(i,y,z),F_3(i,y,z)$. 
Substituting $X,Y,Z$ by $\pm i,y,z$ in the equality~\eqref{eq:y} of Proposition~\ref{prop:degree}, we get $y=\pm 6i/4=\pm 3i/2$. Substituting $X,Y,Z$ by $\pm i,\pm 3i/2,z$ in the equality~\eqref{eq:z} of the same proposition, we get $z=(\pm 2i+3(\pm 3i/2))/\mp i=- 13/2$. Therefore, $\varphi_F^{-1}(\{(1/4,0,0)\})=\{(0,0,1/4),(i,3i/2,-13/2),(-i,-3i/2,-13/2)\}$.

Using {\sc Singular} \cite{dgps} in \url{https://singular-in-browser.pages.dev/}:
\begin{verbatim}
ring A=(0,i),(x,y,z),dp; minpoly =i^2+1;//define Q(i)[x,y,z]
poly f1=(1+x*y)^3*z+y^2*(1+x*y)*(4+3*x*y); //define f1
poly f2=y+3*x*(1+x*y)^2*z+3*x*y^2*(4+3*x*y); //define f2
poly f3=2*x-3*x^2*y-x^3*z; //define f3
ideal I=f1,f2,f3;
subst(I,x,0,y,0,z,1/4); // image of point (0,0,1/4)
subst(I,x,i,y,3i/2,z,-13/2); // image of point (i,3i/2,-13/2)
subst(I,x,-i,y,-3i/2,z,-13/2); // image of point (-i,-3i/2,-13/2)
\end{verbatim}
%\vspace*{-0.7cm}
\end{proof}

\begin{proposition}\label{prop:2max}
Let $F_1,F_2,F_3$ be the three polynomials in {\rm Example~\ref{example:ccj}} and
$F=(F_1,F_2,F_3)$. Let $P=(P_1,P_2,P_3)=(X,Y,Z+\tfrac{1}{4})$, 
$Q=(Q_1,Q_2,Q_3)=(X-\tfrac{1}{4},Y,Z)$ and $G=Q(F(P))$, where
$F,P,Q,G$ are in $\mbk[X,Y,Z]^3$.
Let $\varphi_G:\mbk^3\to \mbk^3$ be the polynomial map defined by $G$. 
Then, $\varphi_G=\varphi_Q\circ\varphi_F\circ\varphi_P$ is a non injective Keller map. 
\begin{itemize}
\item[$(a)$] If $\mbk=\mbc$, then $I_G=\langle X,Y,Z\rangle\cap
\langle X-i,Y-\tfrac{3}{2}i,Z+\tfrac{27}{4}\rangle
\cap \langle X+i,Y+\tfrac{3}{2}i,Z+\tfrac{27}{4}\rangle$. In particular, 
$V_G=\{(0,0,0),(i,\frac{3}{2}i,-\tfrac{27}{4}),(-i,-\tfrac{3}{2}i,-\tfrac{27}{4)}\}$ and
$\card(V_G)=v_G=[\mbc(X,Y,Z):\mbc(G)]=3$.
\item[$(b)$] If $\mbk=\mbr$, then $I_G=\langle X,Y,Z\rangle\cap
\langle X-\tfrac{2}{3}Y,Y^2+\tfrac{9}{4},Z+\tfrac{27}{4}\rangle$. In particular, 
$V_G=\{(0,0,0)\}$ and $1=\card(V_G)<v_G=2<[\mbr(X,Y,Z):\mbr(G)]=3$.
\end{itemize}
\end{proposition}
\begin{proof}
Set $a=(0,0,1/4)$, $b=(i,3i/2,-13/2)$, $b'=(-i,-3i/2,-13/2)$ and $c=(1/4,0,0)$.
By Lemma~\ref{lem:2max}, $\varphi_F(a)=\varphi_F(b)=\varphi_F(b')=c$. 
By Lemma~\ref{lem:notinjective}, we deduce that 
$\varphi_G=\varphi_Q\circ\varphi_F\circ\varphi_P$ is a non injective Keller map, 
$V_G\supseteq \{0,b-a,b'-a\}$ and $I_G\subseteq\mfm_1\cap\mfm_2\cap \mfm_3$, where 
$\mfm_1=\langle X,Y,Z\rangle$, 
$\mfm_2=\langle X-i,Y-\tfrac{3}{2}i,Z+\tfrac{27}{4}\rangle$
and $\mfm_3=\langle X+i,Y+\tfrac{3}{2}i,Z+\tfrac{27}{4}\rangle$.
The equality $I_G=\mfm_1\cap\mfm_2\cap \mfm_3$ follows as in the proof of 
Proposition~\ref{prop:IF<xyz}. Since $X-\tfrac{2}{3}Y=(X\mp i)-\tfrac{2}{3}(Y\mp\tfrac{3}{2}i)$ and $Y^2+\tfrac{9}{4}=(Y-\tfrac{3}{2}i)(Y+\tfrac{3}{2}i)$, we deduce that
$J:=\langle X-\tfrac{2}{3}Y,Y^2+\tfrac{9}{4},Z+\tfrac{27}{4}\rangle\subseteq \mfm_2\cap\mfm_3$.
This inclusion induces a surjective morphism of $\mbk$-algebras:
\begin{eqnarray*}
\mbk[Y]/(Y^2+\tfrac{9}{4})\cong \mbk[X,Y,Z]/J\to\mbk[X,Y,Z]/(\mfm_2\cap\mfm_3)\cong 
(\mbk[X,Y,Z]/\mfm_2)\times(\mbk[X,Y,Z]/\mfm_3).
\end{eqnarray*}
Since both rings have dimension $2$ as $\mbk$-vector spaces, we get the equality 
$J=\mfm_2\cap\mfm_3$. 

To check this we can use
{\sc Singular} \cite{dgps}, \url{https://singular-in-browser.pages.dev/}:
\begin{verbatim}
ring A=(0,i),(x,y,z),dp;  minpoly =i^2+1; //define Q(i)[x,y,z]
option(redSB); // forcing the unicity of the reduced GB
poly f1=(1+x*y)^3*z+y^2*(1+x*y)*(4+3*x*y); //define f1 
poly f2=y+3*x*(1+x*y)^2*z+3*x*y^2*(4+3*x*y); //define f2
poly f3=2*x-3*x^2*y-x^3*z; //define f3
poly f1p=subst(f1,z,z+1/4); //define f_1(p)
poly f2p=subst(f2,z,z+1/4); //define f_2(p)
poly f3p=subst(f3,z,z+1/4); //define f_3(p)
poly g1=f1p-1/4; poly g2=f2p; poly g3=f3p; ideal Ig=g1,g2,g3; //define ideal Ig
ideal m1=x,y,z; ideal m2=x-i,y-(3i/2),z+27/4; ideal m3=x+i,y+(3i/2),z+27/4;
ideal In=intersect(m1,m2,m3); std(Ig); std(In); intersect(m2,m3);
\end{verbatim}
%\vspace*{-0.8cm}
\end{proof}

\section{The rational case}

If $\mbk=\mbq$, there exists a Keller map $\varphi_G$ with
$0=\card(V_G)<v_G=1<[\mbk(X,Y,Z):\mbk(G)]=3$. 

\begin{lemma}\label{lem:casQ}
Let $F_1,F_2,F_3$ be the three polynomials in {\rm Example~\ref{example:ccj}} and $F=(F_1,F_2,F_3)$. 
Let $A,B,C\in\mbk[Y_1,Y_2,Y_3]$ be the polynomials defined in {\rm Proposition~\ref{prop:degree}}.
If $\mbk=\mbc$, then $(\tfrac{1}{3},\tfrac{4}{3},1)$ is in the image of $\varphi_F$.  If $\mbk=\mbq$, then $(\tfrac{1}{3},\tfrac{4}{3},1)$ is not in the image of $\varphi_F$ and $A(\tfrac{1}{3},\tfrac{4}{3},1)=\frac{25}{27}$, $B(\tfrac{1}{3},\tfrac{4}{3},1)=0$ and $C(\tfrac{1}{3},\tfrac{4}{3},1)=-2$.
\end{lemma}
\begin{proof}
If $(a,b,c)=\varphi_F(x,y,z)$, for some $(x,y,z)\in\mbk^3$, then, 
by the equality~\eqref{eq:resultant} in Corollary~\ref{cor:inIG}:
\begin{eqnarray}\label{eq:lemcasQ}
A(a,b,c)x^3+B(a,b,c)x+C(a,b,c)=0,
\end{eqnarray}
where $A(a,b,c)=(27a^2c^2-18abc+16a+b^3c-b^2)$, $B(a,b,c)=4-3bc$ and $C(a,b,c)=-2c$. Let us choose $(a,b,c)$ such that $B(a,b,c)=0$ and $A(a,b,c)X^3+C(a,b,c)$ be irreducible in $\mbq[X]$. So $b=4/3$. Taking $c=1$, then $A(a,b,c)=27a^2-8a+16/27$. If $a=1/3$, then $A(a,b,c)=25/27$ and $\tfrac{25}{27}X^3-2$ is irreducible in $\mbq[X]$. Thus, for $(a,b,c)=(1/3,4/3,1)$, there is no rational $x$ satisfying equation~\eqref{eq:lemcasQ}. It follows that the point $(1/3,4/3,1)$ is not in the image of $\varphi_F$. 
Let $\mbk=\mbc$. Suppose that $(1/3,4/3,1)=\varphi_F(x,y,z)=(F_1(x,y,z),F_2(x,y,z),F_3(x,y,z))$, where
$x^3=54/25$, so $1/x=(25/54)x^2$ and $1/x^3=(25/54)$. Substituting $X,Y,Z$ by $x,y,z$ in the equalities~\eqref{eq:y} and \eqref{eq:z} of Proposition~\ref{prop:degree}, 
we get $y=\tfrac{-25x^2+30x+36}{54}$ and $z=\tfrac{-50x^2+125x-115}{54}$. Therefore, 
$\varphi_F^{-1}(\{(\tfrac{1}{3},\tfrac{4}{3},1)\}))=\{(x_k,\tfrac{-25x_k^2+30x_k+36}{54},\tfrac{-50x_k^2+125x_k-115}{54})\mid x_k=e^{2k\pi i/3}\sqrt[3]{\tfrac{54}{25}}\}$.

With {\sc Singular} \cite{dgps}, \url{https://singular-in-browser.pages.dev/}:
\begin{verbatim}
ring A=(0,a),(x,y,z),dp; 
minpoly =a^3-54/25; // defines Q(a) with a^3=54/25
poly f1=(1+x*y)^3*z+y^2*(1+x*y)*(4+3*x*y); //define f1 
poly f2=y+3*x*(1+x*y)^2*z+3*x*y^2*(4+3*x*y); //define f2
poly f3=2*x-3*x^2*y-x^3*z; //define f3
ideal I=f1,f2,f3; 
poly xk=a; poly yk=(-25*a^2+30*a+36)/54; poly zk=(-50*a^2+125*a-115)/54;
subst(I,x,xk,y,yk,z,zk); //returns 1/3,4/3,1
\end{verbatim}
\end{proof}

\begin{proposition}\label{prop:casQ}
Let $F_1,F_2,F_3$ be the three polynomials in {\rm Example~\ref{example:ccj}} and $F=(F_1,F_2,F_3)$. 
Let $G_1=F_1-\tfrac{1}{3}$, $G_2=F_2-\tfrac{4}{3}$, $G_3=F_3-1$ and $G=(G_1,G_2,G_3)$. 
Let $\varphi_G:\mbk^3\to \mbk^3$ be the polynomial map defined by $G$. Then, $\varphi_G$ is a Keller map. If $\mbk=\mbc$, $I_G$ is proper. If $\mbk=\mbq$, then $V_G$ is empty and $I_G$ is a maximal ideal of $\mbq[X]$. In particular, $0=\card(V_G)<v_G=1<[\mbq(X,Y,Z):\mbq(G)]=3$.
\end{proposition}
\begin{proof} 
By Example~\ref{example:ccj}, $\varphi_F$ is a Keller map. Let $Q=(X-\tfrac{1}{3},Y-\tfrac{4}{3},Z-1)\in\mbk[X,Y,Z]^3$. Then, $G=Q(F)$ and $\varphi_G=\varphi_Q\circ\varphi_F$. Since $\varphi_F$ and $\varphi_Q$ are Keller maps, then $\varphi_G$ is a Keller map. 

Suppose that $\mbk=\mbc$. By Lemma~\ref{lem:casQ}, there exists $(x,y,z)\in\mbc^3$ with $\varphi_F(x,y,z)=(\tfrac{1}{3},\tfrac{4}{3},1)$. Hence, $\varphi_G(x,y,z)=\varphi_Q(\varphi_F(x,y,z))=
\varphi_Q(\tfrac{1}{3},\tfrac{4}{3},1)=(0,0,0)$ and $(x,y,z)\in \varphi_G^{-1}(\{(0,0,0)\})=V_G$. Thus, 
$I(\{(x,y,z)\})\supset I(V_G)\supset I_G$, where $I(\{(x,y,z)\})$ is maximal, and so $I_G$ is proper.
(In fact, $I_G$ is the intersection of three maximal ideals.)

Suppose that $\mbk=\mbq$. By Lemma~\ref{lem:casQ}, $(\tfrac{1}{3},\tfrac{4}{3},1)\not\in\im(\varphi_F)$. 
Then, by Lemma~\ref{lem:vgnobuit}, $V_G$ is empty. 

Since $\langle G_1,G_2,G_3\rangle_{\mathbb{C}[X,Y,Z]}$ is a proper ideal of $\mbc[X,Y,Z]$, then, $I_G=\langle G_1,G_2,G_3\rangle_{\mathbb{Q}[X,Y,Z]}$
is a proper ideal of $\mbq[X,Y,Z]$, otherwise $1=P_1G_1+P_2G_2+P_3G_3$, for some $P_i\in\mbq[X,Y,Z]\subset\mbc[X,Y,Z]$ and $1$ would be in $I_{G}\mbc[X,Y,Z]$.

By Corollary~\ref{cor:inIG} and Lemma~\ref{lem:casQ}, 
$A(\tfrac{1}{3},\tfrac{4}{3},1)X^3+B(\tfrac{1}{3},\tfrac{4}{3},1)X+C(\tfrac{1}{3},\tfrac{4}{3},1)=
\frac{25}{27}X^3-2\in I_G$. Consider the morphism of $\mbq$-algebras $\phi:\mbq[T]\to\mbq[X,Y,Z]$ defined by $\phi(T)=X$. Since $\tfrac{25}{27}X^3-2$ is in $I_G$ and $I_G$ is proper, $\phi$ induces
the morphism of $\mbq$-algebras:
\begin{eqnarray*}
\overline{\phi}:K\stackrel{\rm def}{=}\mbq[T]/\langle\tfrac{25}{27}T^3-2\rangle\to 
\mbq[X,Y,Z]/I_G,
\end{eqnarray*}
where $\overline{\phi}(T)=\overline{X}$. Clearly, $K$ is a field, so $\overline{\phi}$ is injective, and $K$ is a $\mbq$-vector space of dimension 3. In particular, by Theorem~\ref{theo:mpd} and
Propositions~\ref{prop:algext} and \ref{prop:degree}, 
\begin{eqnarray*}
3\leq \dim_{\mathbb{Q}}(\mbq[X,Y,Z]/I_G)\leq [\mbq(X,Y,Z):\mbq(G)]=[\mbq(X,Y,Z):\mbq(F)]=3.
\end{eqnarray*}
Thus, $\dim_{\mathbb{Q}}(\mbq[X,Y,Z]/I_G)=3$. Since $\overline{\phi}$ is an injective linear map between two vector spaces of dimension $3$, it is also surjective. Therefore, $\overline{\phi}$ is an isomorphism of $\mbq$-algebras. Since $K$ is a field, it follows that $\mbq[X,Y,Z]/I_G$ is a field and $I_G$ is a maximal ideal. Thus, $0=\card(V_G)<v_G=1<[\mbq(X,Y,Z):\mbq(G)]=3$.

To check $I_G$ is maximal with {\sc Singular} \cite{dgps}, \url{https://singular-in-browser.pages.dev/}:
\begin{verbatim}
LIB "primdec.lib";
ring A=0,(x,y,z),dp;  //define the ring
poly f1=(1+x*y)^3*z+y^2*(1+x*y)*(4+3*x*y); poly g1=f1-(1/3); //define f1 and g1
poly f2=y+3*x*(1+x*y)^2*z+3*x*y^2*(4+3*x*y); poly g2=f2-(4/3); //define f2 and g2
poly f3=2*x-3*x^2*y-x^3*z; poly g3=f3-1; //define f3 and g3
ideal Ig=g1,g2,g3;  list pd=primdecGTZ(Ig); // primary decomposition of Ig
\end{verbatim}
%\vspace*{-0.5cm}
\end{proof}

\begin{remark}\label{rem:final}
If $\mbk=\mbc$, it is not possible to find an example with $\card(V_G)=v_G=2$, just by
applying some invertible affine transformations on Example~\ref{example:ccj}. This is because of
Theorem~4.2 in \cite{ai}, which states that $\card(\varphi_F^{-1}\{(a,b,c)\})\in\{0,1,3\}$.
For any $\mbk$, it is also not possible to find an example with $\card(V_G)=0<v_G=2$,
just by applying some invertible affine transformations on Example~\ref{example:ccj}. Indeed, 
by Theorem~\ref{theo:mpd}, $2=v_G\leq \dim_{\Bbbk}(\mbk[X,Y,Z]/I_G)\leq [\mbk(X,Y,Z):\mbk(G)]=3$. But, if $v_G=2$, then $I_G=\mfm_1\cap \mfm_2$ and $\mbk[X,Y,Z]/I_G\cong \mbk[X,Y,Z]/\mfm_1\times \mbk[X,Y,Z]/\mfm_2$. Since $\card(V_G)=0$, then $\mfm_i$ is not the maximal ideal corresponding to a point and $\dim_{\Bbbk}(\mbk[X,Y,Z]/\mfm_i)\geq 2$, a contradiction. 
\end{remark}

\section{Topological reformulation of the Jacobian Conjecture}

We finish by giving a purely topological reformulation of the Jacobian Conjecture.
The proof is adapted from the one in \cite[Conjecture C]{mn}, where Miranda-Neto shows that
Conjecture C and the Jacobian Conjecture are equivalent.

\begin{conjecture} 
Let $\mbk$ be a field of characteristic zero. 
Let $\varphi_F:\mbk^n\to\mbk^n$ be a Keller map. Then, $V_F$ is irreducible. 
\end{conjecture}

\begin{remark}
Observe that Proposition~\ref{prop:IFnoPropi} gives a counterexample to Conjecture D when $n=3$.
\end{remark}

\begin{theorem}\label{theo:CD=CJ}
For each $n\geq 1$, Conjecture {\rm D} and the Jacobian Conjecture are equivalent. 
\end{theorem}
\begin{proof}
Suppose that the Jacobian Conjecture is true. Let $\varphi_F$ be a Keller map. 
By hypothesis, $\varphi_F$ is an automorphism of $\mbk^n$. By Lemma~\ref{lem:IFmax}, $V_F$ 
is irreducible. 

Fix $n\geq 1$. Suppose that Conjecture D holds for such $n$. 
Let us show that the Jacobian Conjecture also holds for this same $n$. 
According to \cite[Lemma~1.1.14]{vde1}, which holds (separately) for each specific $n\geq 1$, 
we can suppose that $\mbk=\mbc$. Thus, let $\varphi_F:\mbc^n\to\mbc^n$ be a Keller map, 
where $F=(F_1,\ldots,F_n)\in\mbc[X]^n$. 
Let us see that $\varphi_F$ is invertible. 
By \cite[Proof of Theorem~4.2.1, page 82]{vde1}, it is enough to prove that  
$\varphi_F$ is injective. Take
$a,b\in\mbc^n$ with $\varphi_F(a)=\varphi_F(b)$. 
Let $P=(P_1,\ldots,P_n)=(X_1+a_1,\ldots,X_n+a_n)$,
$Q=(Q_1,\ldots,Q_n)=(X_1-F_1(a),\ldots,X_n-F_n(a))$ and 
$G=Q(F(P))=F(X+a)-F(a)$, where $P,Q,G\in\mbc[X]^n$. By definition, 
$\varphi_G=\varphi_Q\circ\varphi_F\circ\varphi_P$.
Clearly, $\varphi_P,\varphi_Q:\mbc^n\to\mbc^n$ are invertible maps. 
Therefore, $\varphi_G:\mbc^n\to\mbc^n$ is a Keller map. 
Using the hypothesis that the Conjecture D is true, we deduce that $V_G$ is irreducible. 
Since $G(0)=0$, the constant term of each $G_i$ is null, so 
$I_G=\langle G_1,\ldots,G_n\rangle\subseteq \langle X_1,\ldots,X_n\rangle$.
(At this point, here, we can not use Conjecture C since, in this implication, 
we are not assuming that it holds.)
By Miranda-Neto's Theorem~3.1 in \cite{mn}, $I_G=\mbc[X]$ or $I_G=N_1\cap\ldots\cap N_s$ 
is a radical ideal intersection of a finite number of maximal ideals $N_i$. 
Since $V_G=V(I_G)$ is irreducible, then $I_G\subsetneq \mbc[X]$ must be a proper ideal 
and $N_i=I(\{p_i\})$, for some $p_i=(p_{i,1},\ldots,p_{i,n})\in\mbc^n$.
Thus, $I_G=N_1\cap\ldots\cap N_s$ and $V_G=V(N_1\cap\ldots\cap N_s)=
V(N_1)\cup\ldots\cup V(N_s)=V(I(\{p_1\}))\cup\ldots\cup V(\{p_s\})=\{p_1,\ldots,p_s\}$. 
Again, since $V_G$ is irreducible, then $s=1$ and $I_G=N_1$ is a maximal ideal. 
But, since $N_1=I_G\subseteq \langle X_1,\ldots,X_n\rangle$, then
$I_G=\langle X_1,\ldots,X_n\rangle$. This means that, for each $i=1,\ldots, n$, 
$X_i=\sum_{j=1}^ng_{i,j}\cdot G_j$, for some $g_{i,j}\in\mbc[X]$. Evaluating $X$ at $b-a$, 
and using that $G_j(b-a)=F_j((b-a)+a)-F_j(a)=F_j(b)-F_j(a)=0$, one deduces that $b_i-a_i=0$ and so
$a=b$.
\end{proof}

\section*{Acknowledgement}

We are grateful to Jos\'e 
M\textsuperscript{a} Giral, Zbigniew Hajto, Cleto B. Miranda-Neto, Bernat Plans
and Santiago Zarzuela for stimulating conversations which have improved
the readability and the interest of the work. 
Particularly, thanks are due to Zbigniew Hajto who noticed a mistake and several misprints
in a previous version. Computations in this work were performed using the computer algebra system {\sc Singular} \cite{dgps} developed at the University of Kaiserslautern.

The second author is partially supported by project PID2023-146936NB-I00, funded by the Spanish State Agency MICIU/AEI/10.13039/501100011033 and by the Generalitat de Catalunya, Departament de Recerca i Universitats, through AGAUR grant 2025 SGR 00878 (GEOMVAP).

\vspace*{0.2cm}

\noindent {\bf AI Disclosure}.
This work was developed with the assistance of the AI tool \cite{gemini}, which was used to ask for ideas and find simpler approaches to proofs as in Theorem~\ref{theo:mpd}, Propositions~\ref{prop:degree}, \ref{prop:2max} and \ref{prop:casQ} and Lemmas~\ref{lem:vgnobuit} and \ref{lem:notinjective}. We have also used it to improve the lines of code in Singular. To us, this has been a new and unexpected way of working. We take full responsibility for the content of this article.

\vspace*{0.2cm}

\noindent {\bf Conflict of interest}. The authors declare that they have no conflict of interest. The authors have no other relevant financial or non-financial interests to disclose.

\end{document}